\documentclass{article}
\usepackage{amsmath,amssymb,latexsym}

\begin{document}
\title{Theorems about the divisors of numbers contained in the form $paa \pm qbb$\footnote{Delivered to the
St.--Petersburg Academy September 23, 1748.
Originally published as {\em Theoremata circa divisores numerorum in hac forma
$paa \pm qbb$ contentorum},
Commentarii academiae scientiarum Petropolitanae \textbf{14} (1751), 151--181,
and
republished in \emph{Leonhard Euler, Opera Omnia}, Series 1:
Opera mathematica,
Volume 2, Birkh\"auser, 1992. A copy of the original text is available
electronically at the Euler Archive, at http://www.eulerarchive.org. This paper
is E164 in the Enestr\"om index.}}
\author{Leonhard Euler\footnote{Date of translation: June 2, 2006.
Translated from the Latin
by Jordan Bell, 4th year undergraduate in Honours Mathematics, School of Mathematics and Statistics, Carleton University,
Ottawa, Ontario, Canada.
Email: jbell3@connect.carleton.ca.}}
\date{}
\maketitle

In the following theorems, the letters $a$ and $b$ designate arbitrary 
relatively prime integers, that is, which aside from
unity have no other common divisor. 

\begin{center}{\Large Theorem 1.}\end{center}
All the prime divisors of numbers contained in the form $aa+bb$ are either
2 or are numbers of the form $4m+1$.

\begin{center}{\Large Theorem 2.}\end{center}
All prime numbers of the form $4m+1$ in turn are contained in the form
of numbers $aa+bb$.

\begin{center}{\Large Theorem 3.}\end{center}
Thus the sum of two squares, that is a number of the form $aa+bb$, is never
able to be divided by any number of the form $4m-1$.

\begin{center}{\Large Theorem 4.}\end{center}
All the prime divisors of numbers contained in the form $aa+2bb$ are either
2, or numbers contained in the form $8m+1$ or in the form $8m+3$.

\begin{center}{\Large Theorem 5.}\end{center}
All prime numbers contained in the forms $8m+1$ or $8m+3$ in turn are numbers
of the form $aa+2bb$.

\begin{center}{\Large Theorem 6.}\end{center}
No number of the form $aa+2bb$ is able to be divided by a number of the forms
$8m-1$ or $8m-3$.

\begin{center}{\Large Theorem 7.}\end{center}
All the prime divisors of numbers contained in the form $aa+3bb$ are either 2 or
3, or are contained in one of the forms $12+1, 12m+7$.

\begin{center}{\Large Theorem 8.}\end{center}
All prime numbers contained in either of the forms $12m+1$ or $12m+7$,
that is 
in the
one form $6m+1$, are numbers of the form $aa+3bb$.

\begin{center}{\Large Theorem 9.}\end{center}
No number of the form $12m-1$ or of the form $12m-7$, that is no number of the
form $6m-1$, is a divisor of any number contained in the form $aa+3bb$.

\begin{center}{\Large Theorem 10.}\end{center}
All the prime divisors of numbers contained in the form $aa+5bb$ are
either 2, or 5, or are contained in one of 4 forms $2m+1,20m+3,20m+7,20m+9$.

\begin{center}{\Large Theorem 11.}\end{center}
If a number $20m+1,20m+3,20m+9,20m+7$ were prime, then it will follow that
\begin{align*}
20m+1=aa+5bb;&&2(20m+3)=aa+5bb\\
20m+9=aa+5bb;&&2(20m+7)=aa+5bb
\end{align*}

\begin{center}{\Large Theorem 12.}\end{center}
No number contained in one of the following forms $20m-1; 20m-3; 20m-9; 20m-7$
is able to be a divisor of any number of the form $aa+5bb$.

\begin{center}{\Large Theorem 13.}\end{center}
All the prime divisors of numbers contained in the form $aa+7bb$ are either
2 or 7 or are contained in one of the following six forms
\[
\begin{array}{ll|l}
&&\textrm{that is in one of these}\\
28m+1&28m+11&14m+1\\
28m+9&28m+15&14m+9\\
28m+25&28m+23&14m+11
\end{array}
\]

\begin{center}{\Large Theorem 14.}\end{center}
If a prime number were contained in one of the forms $14+1,14m+9,14m+11$, then
at once it would be contained in the form $aa+7bb$.

\begin{center}{\Large Theorem 15.}\end{center}
No number of the form $aa+7bb$ is able to be divided by any number which
is contained in one of the following six forms
\[
\begin{array}{ll|l}
&&\textrm{that is of these three}\\
28m+3,&28m+5&14m+3\\
28m+13,&28m+17&14m+5\\
28m+19,&28m+27&14m+13
\end{array}
\]

\begin{center}{\Large Theorem 16.}\end{center}
All the prime divisors of numbers contained in the form $aa+11bb$ are either
2 or 11, or are contained in one of the following 
\[
\begin{array}{ll|l}
\textrm{10 forms}&&\textrm{5 forms}\\
44m+1&44m+3&22m+1\\
44m+9&44m+27&22m+3\\
44m+37&44m+23&22m+9\\
44m+25&44m+31&22m+5\\
44m+5&44m+15&22m+15
\end{array}
\]

\begin{center}{\Large Theorem 17.}\end{center}
If a prime number were contained in any of these ten or 
five forms, then at once either itself or four times it will be a number of
the form $aa+11bb$.

\begin{center}{\Large Theorem 18.}\end{center}
No number of the form $aa+11bb$ is able to be divided by any number which
is contained in one of the following
\[
\begin{array}{ll|l}
\textrm{10 forms}&&\textrm{5 forms}\\
44m+7&44m+29&22m+7\\
44m+13&44m+35&22m+13\\
44m+17&44m+39&22m+17\\
44m+19&44m+41&22m+19\\
44m+21&44m+42&22m+21
\end{array}
\]

\begin{center}{\Large Theorem 19.}\end{center}
All the prime divisors of numbers contained in the form $aa+13bb$ are
either 2 or 3 or are contained in one of the following 12 formulas.
\begin{align*}
52m+1&&52m+7\\
52m+49&&52m+31\\
52m+9&&52m+11\\
52m+25&&52m+19\\
52m+29&&52m+47\\
52m+17&&52m+15
\end{align*}

\begin{center}{\Large Theorem 20.}\end{center}
All prime numbers which are contained in the first of the above
columns of formulas 
are at once numbers of the form $aa+13bb$. On the other hand, twice the
prime numbers
which are contained in the second column of formulas are numbers of the form
$aa+13bb$.

\begin{center}{\Large Theorem 21.}\end{center}
No number of the form $aa+13bb$ is able to be divided by any number which is
contained in one of the following formulas
\begin{align*}
52m+3&&52m+35\\
52m+5&&52m+37\\
52m+21&&52m+41\\
52m+23&&52m+43\\
52m+27&&52m+45\\
52m+33&&52m+51
\end{align*}

\begin{center}{\Large Theorem 22.}\end{center}
All the prime divisors of numbers contained in the form $aa+17bb$ are either
2 or 17 or are contained in one of the following forms
\begin{align*}
68m+1&&68m+3\\
68m+9&&68m+27\\
68m+13&&68m+39\\
68m+49&&68m+11\\
68m+33&&8m+31\\
68m+25&&68m+7\\
68m+21&&68m+63\\
68m+53&&68m+23
\end{align*}

\begin{center}{\Large Theorem 23.}\end{center}
All the prime numbers which are contained in the first
of the above columns of formulas are either themselves or nine times them
numbers of the form $aa+17bb$.
On the other hand, three times the prime numbers in the other column are numbers
of the form $aa+17bb$.

\begin{center}{\Large Theorem 24.}\end{center}
No number of the form $aa+17bb$ is able to be divided by any number which is
contained
in one of the following formulas
\begin{align*}
68m-1&&68m-3\\
68m-9&&68m-27\\
68m-13&&68m-39\\
68m-49&&68m-11\\
68m-33&&68m-31\\
68m-25&&68m-7\\
68m-21&&68m-63\\
68m-53&&68m-23
\end{align*}

\begin{center}{\Large Theorem 25.}\end{center}
All the prime divisors of numbers contained in the form $aa+19bb$ are either
2, or 19, or are contained in one of the following
\[
\begin{array}{ll|l}
\textrm{18 formulas}&&\textrm{9 formulas}\\
76m+1&76m+5&38m+1\\
76m+25&76m+49&38m+5\\
76m+17&76m+9&38m+7\\
76m+45&76m+73&38m+9\\
76m+61&76m+7&38m+11\\
76m+35&76m+23&38m+17\\
76m+39&76m+43&38m+23\\
76m+63&76m+11&38m+25\\
76m+55&76m+47&38m+35
\end{array}
\]

\begin{center}{\Large Theorem 26.}\end{center}
All prime numbers which are contained in one of these forms are either
themselves or four times them numbers of the form $aa+19bb$.

\begin{center}{\Large Theorem 27.}\end{center}
No number of the form $aa+19bb$ is able to be divided by any number which
may be contained in any of the following 9 formulas
\begin{align*}
38m-1\\
38m-5\\
38m-7\\
38m-9\\
38m-11\\
38m-17\\
38m-23\\
38m-25\\
38m-35
\end{align*}
Thus the character of the forms $aa+qbb$ is contained in these theorems;
if $q$ is a prime number, we see first for
all the prime divisors of these forms to be 2 or $q$, or to be
able to be expressed in the form $4qm+\alpha$, so that
no divisor is not contained in this form, and indeed also that each prime number
$4qm+\alpha$ is at once a divisor of such a  form $aa+qbb$. Then moreover
it 
can be deduced that that if a prime number of the form $4qm+\alpha$ were
 a divisor of some number $aa+qbb$, then no number of the form $4qm-\alpha$ 
 can be a divisor of the expression $aa+qbb$. It is clear
 from this therefore that as $4mq+1$ is always contained among the divisors
 of the form $aa+qbb$,
no number $aa+qbb$ is able to be divided
 by any number of the form $4mq-1$. Indeed it is clear with attention
 that were $q$ a prime number of the form $4n-1$ for the forms 
 of the divisors to be able to be reduced to less than twice the number,
so that they can be reduced
 to the formulas $2qm+\alpha$, and that this cannot be done if
 $q$ were a prime number of the form $4n+1$. If therefore
 $4(4n+1)m+\alpha$ were a divisor for the form $aa+(4n+1)bb$, then
 no number of the form $4(4n+1)m+2(4n+1)+\alpha$ will be able
 to be a divisor of the expression $aa+(4n+1)bb$. We will presently
 give a number of notes about the forms $aa+qbb$, in which
 we shall contemplate when $q$ is not a prime number.

\begin{center}{\Large Theorem 28.}\end{center}
All the prime divisors of numbers contained in the form $aa+6bb$ or in
the form $2aa+3bb$ are either 2 or 3 or are contained in one of the following
formulas
\begin{align*}
24m+1&&24m+7\\
24m+5&&24m+11
\end{align*}

\begin{center}{\Large Theorem 29.}\end{center}
All prime numbers either of the form $24m+1$ or $24m+7$ are contained
in the expression $aa+6bb$; whereas the prime numbers of the form $24m+5$ and
$24m+11$ are contained in the expression $2aa+3bb$.

\begin{center}{\Large Theorem 30.}\end{center}
No number $aa+6bb$ or $2aa+3bb$ is able to be divided by any number which
is contained in any of the following forms
\begin{align*}
24m-1&&24m-5\\
24m-7&&24m-11
\end{align*}

\begin{center}{\Large Theorem 31.}\end{center}
All the prime divisors of numbers in the form $aa+10bb$ or in the form
$2aa+5bb$ are either 2 or 5 or are contained in one of the following
forms
\begin{align*}
40m+1&&40m+7\\
40m+9&&40m+23\\
40m+11&&40m+37\\
40m+19&&40m+13
\end{align*}

\begin{center}{\Large Theorem 32.}\end{center}
Prime numbers contained in the forms of the first column above are 
at once numbers of the form $aa+10bb$ and prime numbers contained in the second
column are
numbers of the form $2aa+5bb$.

\begin{center}{\Large Theorem 33.}\end{center}
No number either of the form $aa+10bb$ or of the form $2aa+5$ is able to be
divided by any number which is contained in the following forms
\begin{align*}
40m-1&&40m-7\\
40m-9&&40m-23\\
40m-11&&40m-37\\
40m-19&&40m-13
\end{align*}

\begin{center}{\Large Theorem 34.}\end{center}
All the prime divisors of numbers contained in the form $aa+14bb$ or in
the form $2aa+7bb$ are either 2 or 7 or are contained in one
of the following formulas
\begin{align*}
56m+1&&56m+3\\
56m+9&&56m+27\\
56m+25&&56m+19\\
56m+15&&25m+5\\
56m+23&&56m+45\\
56m+39&&56m+13
\end{align*}

\begin{center}{\Large Theorem 35.}\end{center}
Prime numbers of the formulas contained in the first column above are at once
numbers of either the form $aa+14bb$ or $2aa+7bb$, while on the other hand
those three times those contained in the second column is contained
in the formulas in the first column. 

\begin{center}{\Large Theorem 36.}\end{center}
If in the above examples the sign $+$ is switched to $-$, then no number
contained in this form can be a divisor of either the form $aa+14bb$
or the form $2aa+7bb$.

\begin{center}{\Large Theorem 37.}\end{center}
All the prime divisors of numbers contained in the form $aa+15bb$ or in the form
$3aa+5bb$ are either 2, or 3, or 5 or are contained in one
of the following formulas
\[
\begin{array}{ll|l}
&&\textrm{or of these 4}\\
60m+1&60m+31&30m+1\\
60m+17&60m+47&30m+17\\
60m+19&60m+49&30m+19\\
60m+23&60m+53&30m+23
\end{array}
\]

\begin{center}{\Large Theorem 38.}\end{center}
All prime divisors of numbers contained in the form
$aa+21bb$ or in the form $3aa+7bb$ are either 2, or 3, or 7
or are contained in one of the following formulas
\begin{align*}
84m+1&&84m+5\\
84m+25&&84m+41\\
84m+37&&84m+17\\
84m+55&&84m+11\\
84m+31&&84m+13\\
84m+19&&84m+71
\end{align*}

\begin{center}{\Large Theorem 39.}\end{center}
All the prime divisors of numbers contained in the form $aa+35bb$ or in the form
$5aa+7bb$ 
are either 2, or 5, or 7 or are contained in one of the following formulas
\[
\begin{array}{ll|l}
&&\textrm{or of these}\\
140m+1&140m+3&70m+1\\
140m+9&140m+27&70m+3\\
140m+81&140m+103&70m+9\\
140m+29&140m+87&70m+11\\
140m+121&140m+83&70m+13\\
140m+109&140m+47&70m+17\\
140m+11&140m+33&70m+27\\
140m+99&140m+17&70m+29\\
140m+51&140m+13&70m+33\\
140m+39&140m+117&70m+39\\
140m+71&140m+73&70m+47\\
140m+79&140m+97&70m+51
\end{array}
\]

\begin{center}{\Large Theorem 40.}\end{center}
All the prime divisors of numbers contained in any of the forms
\begin{eqnarray*}
aa+30bb;&&2aa+15bb;\\
3aa+10bb;&&5aa+6bb
\end{eqnarray*}
are either 2, or 3, or 5, or are contained in one
of the following formulas.
\begin{eqnarray*}
120m+1;&&120m+11\\
120m+13;&&120m+23\\
120m+49;&&120m+59\\
120m+37;&&120m+47\\
120m+17;&&120m+67\\
120m+101;&&120m+31\\
120m+113;&&120m+43\\
120m+29;&&120m+79
\end{eqnarray*}
These theorems will suffice to formulate the following notes, from which the
nature
of the divisors of formulas of the type $paa+qbb$ will be examined thoroughly.

\begin{center}{\Large Note 1.}\end{center}
The form $paa+qbb$ has no divisors but that which is at once a divisor of
$aa+pqbb$. The rationale for this is clear; for were a number a divisor
of the form $paa+qbb$, the same then divides the form $ppaa+pqbb$, which
is $aa+pqbb$, by putting $a$ in place of $pa$. Thus it suffices
to consider the single form $aa+Nbb$; from this rationale the divisors for
$paa+qbb$ are completed.

\begin{center}{\Large Note 2.}\end{center}
Among the prime numbers which divide any of the numbers contained in the form
$aa+Nbb$,
the prime 2 occurs.
For if $N$ were an odd number
with $a$ and $b$ taken as odd numbers, the form $aa+Nbb$ will be
divisible by 2; and if $N$ were an even number, with $a$ taken
to be even, the form again will be divisible by 2. Then indeed
this number $N$ will be able to be a divisor of the form $aa+Nbb$, which
by taking $a=N$ is clear.

\begin{center}{\Large Note 3.}\end{center}
All the remaining prime divisors of the form $aa+Nbb$ are thus
able to be expressed as $4Nm+\alpha$; moreover in turn, all prime numbers
contained in the form $4Nm+\alpha$ are at once divisors of the
form $aa+Nbb$. In addition, if the expression $4Nm+\alpha$ permits
divisors of the form $aa+Nbb$, then no number of the form $4Nm-\alpha$ will
be able to be a divisor of any number contained in the form $aa+Nbb$.

\begin{center}{\Large Note 4.}\end{center}
It will be moreover that the particular values of $\alpha$ depend on the
character
of the number $N$; and indeed always unity will be among the
values for $\alpha$. Then also, because prime numbers
are being sought for in the formula $4Nm+\alpha$, it is clear that
no even number nor any number which has a common divisor with $N$ can 
constitute a value of $\alpha$.

\begin{center}{\Large Note 5.}\end{center}
As well, all the values of $\alpha$ will be less than $4N$, for if they were
greater, by decreasing the number $m$, ones less than $4N$ can be obtained.
Hence the values of $\alpha$ will be odd numbers less than $4N$, and
also prime to $N$. But indeed not all of these odd numbers prime to $N$ will
furnish suitable values for $\alpha$, but half of them are excluded
by the rule that if $x$ were a value of $\alpha$, then $-x$, that is
$4N-x$, may not be a value of it; and in turn if $x$ were not a value
of $\alpha$, then $4N-x$ will be certain to be a value of it.

\begin{center}{\Large Note 6.}\end{center}
So that $4Nm+\alpha$ will contain all the prime divisors of the formula
$aa+Nbb$, a value of this $\alpha$ will be defined in the following way.
Were $p,q,r,s$, etc. distinct prime numbers, excepting 2, which will be 
considered separately, then
\[
\begin{array}{lc|l}
\textrm{if it were}&&\textrm{a value of $\alpha$ will be}\\
N=1&&1\\
N=2&&2\\
N=p&&p-1\\
N=2p&&p-1\\
N=2p&&2(p-1)\\
N=pq&&(p-1)(q-1)\\
N=2pq&&2(p-1)(q-1)\\
N=pqr&&(p-1)(q-1)(r-1)\\
N=2pqr&&2(p-1)(q-1)(r-1)\\
&\textrm{etc.}&
\end{array}
\]

\begin{center}{\Large Note 7.}\end{center}
Moreover, in the same way that unity always appears among the values
of $\alpha$, thus too any number which is an odd square and relatively
prime to $N$ shall have a place among the values of $\alpha$.
For by putting $b$ as the even number
$2c$, the formula would be $aa+4Ncc$, which if it were a prime number
must be contained in the expression $4Nm+\alpha$. Therefore $\alpha$ will
be a residue of $aa$, because $aa$ remains when this is divided by $4N$.
In a similar way, among the values of $\alpha$ all numbers of
the form $aa+N$ should appear, that is which remain from division by $4N$;
for by putting $b=2c+1$ it would be $aa+Nbb=aa+N+4N(cc+c)$, which, if it were
a prime number, gives that $aa+N$ is a value of $\alpha$.

\begin{center}{\Large Note 8.}\end{center}
It can be understood as well that if $x$ were a value of $\alpha$, then
too $xx$ (which indeed is clear from the preceeding) and all higher powers
of $x$, so that $x^\mu$ itself should have a place among the values of $\alpha$.
Then, if aside from $x$ also $y$ were a value of $\alpha$, then too
$xy$ and in general $x^\mu y^\nu$ gives a value of $\alpha$. Certainly
if $x^\mu y^\nu$ were greater 
than $4N$, by dividing this the remainder will be a value of $\alpha$.
In a similar way, if in addition $z$ were a value of $\alpha$, then
further $x^\mu y^\nu z^\xi$ will be a value of $\alpha$. And then
from this inquiry, from one or several values of $\alpha$ all the other values
can be found by easy work.

\begin{center}{\Large Note 9.}\end{center}
Were $x$ some number prime to $4N$ and less than it, then either
$+x$ or $-x$ will be a value of $\alpha$. If therefore
$x$ were a prime number, from the following table it can be understood 
whether the case $+x$ or $-x$ for the value of $\alpha$ obtains.
\[
\begin{array}{l|l}
\textrm{If}&\textrm{it will be}\\
N=3n-1&\alpha=+3\\
N=3n+1&\alpha=-3\\
\hline\\
N=\begin{cases}5n+1\\
5n+4\end{cases}&\alpha=+5\\
N=\begin{cases}5n+2\\
5n+3\end{cases}&\alpha=-5\\
\hline\\
N=\begin{cases}7n+3\\7n+5\\7n+6\end{cases}&\alpha=+7\\
N=\begin{cases}7n+1\\7n+2\\7n+4\end{cases}&\alpha=-7\\
\hline\\
N=\begin{cases}11n+2\\11n+6\\11n+7\\11n+8\\11n+10\end{cases}&\alpha=+11\\
N=\begin{cases}11+1\\11n+3\\11n+4\\11n+5\\11n+9\end{cases}&\alpha=-11
\end{array}
\]
If a prime number is given, whether the sign $+$ or $-$ obtains for the value
of $\alpha$ will thus be investigated. Both cases shall be pursued, the one
in which the given prime number is of the form $4u+1$, the other in which
the it is of the form $4u-1$. In the first case it will be $\alpha=+(4u+1)$
if it were $N=(4u+1)n+tt$, or $\alpha=-(4u+1)$ if it were
$N \neq (4u+1)n+tt)$. In the latter case however it will be $\alpha=+(4u-1)$
if it were $N \neq (4u-1)n+tt)$ or $\alpha=-(4u-1)$ if $N=(4u-1)n+tt$.
Here it should be noted that in the way as the sign $=$ denotes equality,
so the sign $\neq$ denotes the impossibility of equality. But if
it were moreover in both cases $N=(4u \pm 1)n+s$, it will also
be $N=(4u \pm 1)n+s^\nu$, where
$\nu$ denotes a certain integral number, from which such a table for any prime
number may be constructed without effort.

\begin{center}{\Large Note 10.}\end{center}
Since $4Nm+1$ is among the forms of the prime divisors of $aa+Nbb$, the
expression $aa+Nbb$ will not be able to be divided by any number which is
contained
in the form $4Nm-1$. In a similar way should $4Nm+tt$ exhibit the form of
divisors
of the expression $aa+Nbb$, it follows that no number of the type
$4Nm-tt$ will be able to be a divisor of any number contained in the form
$aa+Nbb$, whenever $a$ and $b$ are relatively
prime numbers. Then from this impossibility will follow this equation
$(4Nm-tt)u=aa+Nbb$, and thus it will be $4Nmu-ttu-Nbb=aa$, 
if indeed the numbers $4Nmu-ttu$ and $Nbb$ were relatively prime,
because with certainty it follows, if $b=1$ and $t=1$, that this obtains.

\begin{center}{\Large Corollary.}\end{center}
No number contained in the form $4abc-b-c$ is ever able to be a square.

\begin{center}{\Large Note 11.}\end{center}
If $N$ were a number of the form $4n-1$, then the forms of the divisors
are reduced to less than twice the number, so that they may be comprised
in the form $2Nm+\alpha$. Indeed if $4Nm+\alpha$ were a form of divisors,
then too $4Nm+2N+\alpha$ will be a form of divisors. In this way were
$2Nm+tt$ a form of divisors, it follows that no number $4Nm-tt$ can be a divisor
of the form $aa+Nbb$. Thus it will be $(2Nm-tt)u=aa+Nbb$; from this development
$N=4n-1$ arises.

\begin{center}{\Large Corollary.}\end{center}
No number of the form $2abc-b-c$, if either $b$ or $c$ are odd numbers of
the form $4n-1$, is ever able to be a square.

\begin{center}{\Large Note 12.}\end{center}
If $N$ we an odd number of the form $4n+1$ or also an oddly even number, then
the forms of the divisors will
not be able to be reduced to less than twice the number. Indeed if $4Nm+\alpha$
were a divisor of the form $aa+Nbb$ then $4Nm+2N+\alpha$ will not be able to be
a divisor of
the same form. Thus $2(2m+1)N+tt$ will not be a divisor of the form $aa-Nbb$,
and thus the equation $(2(2m+1)N+tt)u=aa+Nbb$ will be impossible if indeed
$a$ and $b$ were relatively prime numbers: and $N$ were an odd number of the
form $4n+1$ or an oddly even number. From which follows this

\begin{center}{\Large Corollary.}\end{center}
No number of the form $2abc-b+c$, with $a$ arising as an odd number and
$b$ either oddly even or odd of the form $4n+1$, is ever able to be a square.

\begin{center}{\Large Scholion 1.}\end{center}
These show that the character of the divisors of the formulas $aa+Nbb$
has been developed sufficiently, and at once take care of
all the forms of divisors which are to be found without burden,
by which having noted the forms of the numbers become known
which are never able to present divisors for the formula
$aa+Nbb$. Therefore by this all the values of $N$ may be known, either
prime or composite; it remains however that we pursue the case in which
$N$ denotes a negative number, either prime or composite; indeed it is
clear for the formula $paa-qbb$ to have no divisor which is not a divisor of
$aa-pqbb$, that is $pqaa-bb$, from which it suffices to pursue just  the
formas $aa-Nbb$.  

\begin{center}{\Large Theorem 41.}\end{center}
All the prime divisors of numbers contained in the form $aa-bb$ are either 2 or
$4m\pm 1$, namely no number is permitted which is not a divisor of two
different squares. In turn moreover all numbers except those which are oddly even
are themselves the difference of two squares.

\begin{center}{\Large Theorem 42.}\end{center}
All the prime divisors of numbers contained in the form $aa-2bb$ are either
2 or of the form $8m \pm 1$. And all prime numbers of the form $8m\pm 1$
are contained infinitely many ways in the form $aa-2bb$.

\begin{center}{\Large Theorem 43.}\end{center}
All the prime divisors of numbers contained in the from $aa-3bb$ are either 2
or 3 or are of the form $12m \pm 1$. And also in turn all prime numbers of
this type are contained at once in the form $aa-3bb$ or in the form $3aa-bb$
infinitely many ways.

\begin{center}{\Large Theorem 44.}\end{center}
All the prime divisors of the form $aa-5bb$ are either 2 or 5 or are contained
\[
\begin{array}{ll|l}
\textrm{in one or the other of these forms}&&\textrm{or in such a one}\\
20m \pm 1&20m\pm 9&10m \pm 1.
\end{array}
\]
And all prime numbers contained in these forms are at once divisors of
the form $aa-5bb$.

\begin{center}{\Large Theorem 45.}\end{center}
All the prime divisors of the form $aa-7bb$ are either 2 or 7 or are contained
in one of the following forms:
\[
28m\pm 1; \quad 28m\pm 3; \quad 28m\pm 9
\]
and also in turn all prime numbers contained in these forms are divisors of
the form $aa-7bb$.

\begin{center}{\Large Theorem 46.}\end{center}
All the prime divisors of the form $aa-11bb$ are either 2 or 11 or are contained
in one of the following forms
\[
44m \pm 1; \quad 44\pm 5; \quad 44m\pm 7;\quad 44m\pm 9; \quad 44m\pm 19
\]
and also in turn all prime numbers contained in these forms are at once
divisors of the form $aa-11bb$; this reciprocity holds in all the following
theorems.

\begin{center}{\Large Theorem 47.}\end{center}
All the prime divisors of the form $aa-13bb$ are either 2 or 13 or are contained
in
one of the following forms:
\[
\begin{array}{ll|l}
&&\textrm{which are reduced to these}\\
52m\pm1&52m\pm 3&26m\pm 1\\
52m\pm 9&52m\pm 25&26m\pm 3\\
52m\pm 23&52m\pm 17&26\pm 9.
\end{array}
\]

\begin{center}{\Large Theorem 48.}\end{center}
All the prime divisors of numbers of the form $aa-17bb$ are either 2 or 17,
or are contained in one of the following forms:
\[
\begin{array}{ll|l}
&&\textrm{which are reduced to these}\\
68m\pm 1&68m\pm 9&34m\pm 1\\
68m\pm 13&68m\pm 19&34m\pm 9\\
68m\pm 33&68m\pm 25&34m\pm 13\\
68m\pm 21&68m\pm 15&34m\pm 15.
\end{array}
\]

\begin{center}{\Large Theorem 49.}\end{center}
All the prime divisors of numbers of the form $aa-19bb$ are either 2 or 19 or
are contained in one of the following forms
\[
\begin{array}{lll}
76m\pm 1&76m\pm 3&76m\pm 9\\
76m\pm 27&76m\pm 5&76m\pm 15\\
76m\pm 31&76m\pm 17&76m\pm 25.
\end{array}
\]

\begin{center}{\Large Theorem 50.}\end{center}
All the prime divisors of numbers of the form $aa-6bb$ are either 2 or 3 or
are contained in one of these forms:
\[
24m\pm 1; \quad 24\pm 5.
\]

\begin{center}{\Large Theorem 51.}\end{center}
All the prime divisors of numbers of the form $aa-10bb$ are either 2 or 3
or are contained in these forms:
\[
\begin{array}{ll}
40m\pm 1&40m\pm 3\\
40m\pm 9&40m\pm 13.
\end{array}
\]

\begin{center}{\Large Theorem 52.}\end{center}
All the prime divisors of numbers of the form $aa-14bb$ are either 2 or 7 or
are contained in these forms
\[
\begin{array}{lll}
56m\pm 1&56m\pm 5&56m\pm 25\\
56m\pm 13&56m\pm 9&56m\pm 11
\end{array}
\]

\begin{center}{\Large Theorem 53.}\end{center}
All the prime divisors of numbers of the form $aa-22bb$ are either 2 or 11 or are contained in these forms
\[
\begin{array}{lll}
88m \pm 1&88m\pm 3&88m\pm 9\\
88m\pm 27&88m\pm 7&88m\pm 21\\
88m\pm 25&88m\pm 13&88m\pm 39\\
&88m\pm 29.&
\end{array}
\]

\begin{center}{\Large Theorem 54.}\end{center}
All the prime divisors of numbers of the form $aa-15bb$ are either
2 or 3 or 5 or are contained in these forms:
\[
60m\pm 1; \quad 60m\pm 7; \quad 60m\pm 11;\quad 60m\pm 17.
\]

\begin{center}{\Large Theorem 55.}\end{center}
All the prime divisors of numbers of the form $aa-21bb$ are either 2 or
3 or 7 or are contained in these forms:
\[
\begin{array}{ll|l}
&&\textrm{which are reduced to these}\\
84m\pm 1&84m\pm 5&42m\pm 1\\
84m\pm 25&84m \pm 41&42m \pm 5\\
84m\pm 37&84m\pm 17&42m\pm 17
\end{array}
\]

\begin{center}{\Large Theorem 56.}\end{center}
All the prime divisors of numbers of the form $aa-33bb$ are either 2 or 3, or
11 or are contained in these forms
\[
\begin{array}{ll|l}
&&\textrm{which are reduced to these}\\
132m\pm 1&132m\pm 17&66m\pm 1\\
132m\pm 25&132m\pm 29&66m\pm 17\\
132m\pm 35&132m\pm 65&66m\pm 25\\
132m\pm 49&132m\pm 41&66m\pm 29\\
132m\pm 37&132m\pm 31&66m\pm 31
\end{array}
\]

\begin{center}{\Large Theorem 57.}\end{center}
All the prime 
divisors of numbers of the form $aa-35bb$ are either 2 or 5 or 7 or are contained
in these forms:
\[
\begin{array}{lll}
140m\pm 1&140m\pm 9&140m\pm 59\\
140m\pm 29&140m\pm 19&140m\pm 31\\
140m\pm 13&140m\pm 23&140m\pm 67\\
140m\pm 43&140m\pm 33&140m\pm 17
\end{array}
\]

\begin{center}{\Large Theorem 58.}\end{center}
All the prime divisors of numbers of the form $aa-30bb$ are either 2 or 3 or
5 or are contained in these forms
\[
\begin{array}{lll}
120m\pm 1&120m\pm 13&120m\pm 49\\
120m\pm 37&120m\pm 7&120m\pm 29\\
120m\pm 17&120m\pm 19&
\end{array}
\]

\begin{center}{\Large Theorem 59.}\end{center}
All the prime divisors of numbers of the form $aa-105bb$ are either 2 or 3
or 5 or 7 or are contained in these forms
\[
\begin{array}{ll|l}
&&\textrm{which are reduced to these}\\
420m\pm 1&420m\pm 13&210m\pm 1\\
420m\pm 169&420m\pm 97&210m\pm 13\\
420m\pm 23&420m\pm 121&210m\pm 23\\
420m\pm 107&420m\pm 131&210m\pm 41\\
420m\pm 109&420m\pm 157&210m\pm 53\\
420m\pm 59&420m\pm 73&210m\pm 59\\
420m\pm 101&420m\pm 53&210m\pm 73\\
420m\pm 151&420m\pm 137&210m\pm 79\\
420m\pm 89&420m\pm 103&210m\pm 89\\
420m\pm 79&420m\pm 187&210m\pm 97\\
420m\pm 41&420m\pm 113&210m\pm 101\\
420m\pm 209&420m\pm 197&210m\pm 103
\end{array}
\]

\begin{center}{\Large Note 13.}\end{center}
Therefore all the prime divisors of numbers contained in the form $aa-Nbb$
are either 2, are divisors of the number $N$, or are comprised in
the form $4Nm\pm \alpha$. And if $4Nm+\alpha$ were a form of divisors,
then also $4Nm-\alpha$ would be a form of divisors: on the other hand,
for the forms $aa+Nbb$, if $4Nm+\alpha$ were the form of a divisor
then this form does not permit a divisor $4Nm-\alpha$.

\begin{center}{\Large Note 14.}\end{center}
Therefore by taking $4Nm \pm \alpha$ as the general form of divisors of 
numbers contained in the expression $aa-Nbb$, the letter $\alpha$
signifies many different numbers; unity is always contained among these,
then indeed from the discussion on prime divisors, among the
values for $\alpha$ there will be no even number and no divisor of
the number $N$. Then moreover it is clear that all the values
of $\alpha$  can be arranged such that each is less than $2N$. 
For if $4Nm+2N+b$ were a divisor, then by putting $m-1$ in place
of $m$, $4Nm-(2N-b)$ will be a divisor. Therefore the values
of $\alpha$ will be odd numbers prime to $N$, less than $2N$; of
all the odd numbers prime to $N$ and less than $2N$, exactly half
will give suitable values for $\alpha$, while the remaining
exhibit formulas in which no divisor whatsoever is contained. 
Indeed always just as many formulas of divisors will be obtained as those which
are not permissible, except for the single case in which $N=1$.

\begin{center}{\Large Note 15.}\end{center}
For the number of values of $\alpha$ attained for the formula
of the divisors of $4Nm \pm \alpha$, since from the varying of the sign
each formula is doupled, the same rule prevails which I gave above
in note 6.
For in the last theorems, in which it was $N=105,3,5,7$, the number
of values of $\alpha$ will be equal to $2,4,6=48$, or as each formula
would be a pair, the number of formulas would be 24, just as many 
indeed as we have exhibited.

\begin{center}{\Large Note 16.}\end{center}
Moreover just as uinity always appears among the values of $\alpha$, so to
any square numbers which is prime with $4N$ provides a suitable value
for $\alpha$. For by putting $b=2c$, the formulae $aa-Nbb$ becomes
$aa-4Ncc$ or $4Ncc-aa$, from which it is clear for any square number
$aa$, which is prime to $4N$, to exhibit a suitable value for $\alpha$, namely
by taking the residue which remains from the division of $aa$ by $4N$.
In a similar way by putting $b=2c+1$, the formula
$Nbb-aa$ becomes $4N(cc+c)+N-aa$, from which too all the numbers
$N-aa$ or $aa-N$, which are prime to $4N$, provide suitable values
for $\alpha$. Then too it is to be noted that if $x,y,z$ are values of
$\alpha$, then too for $x^\mu y^\nu z^\xi$, and any such products which result
from numbers
$x,y,z$
and any powers, to be able to exhibit values of $\alpha$;
from which knowing one or several values of $\alpha$, all can be
obtained with easy work.

\begin{center}{\Large Note 17.}\end{center}
So that the way in which the values for the letter
$\alpha$ are perpetually obtained may be made clear, the following table
is considered, similarly as that in note 9 is had.
\[
\begin{array}{l|l}
\textrm{It will be}&\textrm{if it were}\\
\alpha=3&N=3n+1\\
\alpha \neq 3&N=3n-1\\
\hline\\
\alpha=5&N=5n\begin{cases}+1\\-1\end{cases}\\
\alpha\neq 5&N=5n\begin{cases}+2\\-1\end{cases}\\
\hline\\
\alpha=7&N=7n\begin{cases}+1\\+2\\-3\end{cases}\\
\alpha\neq 7&N=7n\begin{cases}-1\\-2\\+3\end{cases}\\
\hline\\
\alpha=11&N=11n\begin{cases}+1\\-2\\+3\\+4\\+5\end{cases}\\
\alpha \neq 11&N=11n\begin{cases}-1\\+2\\-3\\-4\\-5\end{cases}\\
\hline\\
\alpha=13&N=13n\begin{cases}+1\\-1\\+3\\-3\\+4\\-4\end{cases}\\
\alpha \neq 13&N=13n\begin{cases}+2\\-1\\+5\\-5\\+6\\-6\end{cases}
\end{array}
\]

\begin{center}{\Large Note 18.}\end{center}
From this table therefore, prime numbers, which provide suitable values
for $\alpha$, are easy able discerned, and inappropriate ones can be rejected.
Namely for a given prime number $p$, all the square numbers in the form
$pn+\theta$ are able to be comprehended, which follow by taking square numbers
for $\theta$, that is residues, which remain after the division of squares by
$p$. Thereby if $N$ were a number of the type $pn+tt$, then among
the forms of the divisors
of $4Nm \pm \alpha$ are the formulas $aa-Nbb$, that is $Nbb-aa$, and it will
be had that $\alpha=p$, but if on the other hand the number $N$ were not
contained
in the form $pn+tt$,
then no number contained in the form $4Nm\pm p$ will be able to be a divisor
of any number of the form $aa-Nbb$.

\begin{center}{\Large Note 19.}\end{center}
If $N$ is an odd number of the form $4n+1$ then the forms of the divisors
of the expression $aa-Nbb$ are able to be reduced to 
twice as small as $4Nm \pm \alpha$, so that they can be presented in
this way: $2Nm\pm \alpha$. Indeed in this case, if $4Nm \pm \alpha$ were
a form of divisors, then too $4Nm \pm (2N-\alpha)$ will be a form
of divisors, in such a way as with the case $N=13$, where $52m \pm 3$
was a form of divisors of $aa-13bb$, and likewise $52m \pm 23$ will be a form
of divisors.

\begin{center}{\Large Note 20.}\end{center}
But if moreover $N$ were an oddly even number, or an odd number
of the form $4n-1$, then this reduction of dividing the form into two smaller
ones does not succeed. Indeed in this case of the formula $aa-Nbb$ if
$4Nm \pm \alpha$ were the form of divisors, then $4Nm \pm (2N-\alpha)$
will not be such, that is: no number contained in the form $2(2m \pm 1)N \pm \alpha$ will be a divisor of any number of the form $aa-Nbb$. Therefore by putting
$\alpha=tt$, it will be
\[
(2(2m \pm 1)N \pm tt) u \neq aa-Nbb,
\]
of which the following is a consequence.

\begin{center}{\Large Corollary.}\end{center}
No number contained in the form $2abc\pm c+b$ is able to be a square, if
$a$ were an odd number, and $b$ were either oddly even number or odd
of the form $4n-1$.

\begin{center}{\Large Scholion 2.}\end{center}
Innumerable special formulas are able to be deduced from the above, which are
never able to be made into squares.
We considered the first form $aa+Nbb$ and it would be a formula of the
type $4Nm+A$,
so that 
no number contained in this form is able to be a divisor of the form $aa+Nbb$.
It will therefore be $aa+Nbb=(4Nm+A)u$,
where the sign $\neq$ denotes that the equation is impossible, from which
it follows that $aa \neq 4Nmu+Au-Nbb$. Were it $b=Ac$ it would be $aa \neq
4Nmu+Au-NAAcc$. It is put in turn $u=NAcc+d$, and it will
be $aa \neq 4NNAmcc+4Nmd+Ad$. Were it $d=4NNn$ it will be
$aa \neq 16N^3 mn+4NNAmcc+4NNAn$. Were this formula
divided by the square $4NN$ and it put $c=1$, it will bethe formula
$4Nmn+Am+An$, which is never able to be a square if the form $aa+Nbb$ is not
able to be divided by any number contained in the form $4Nm+A$. From the above
theorems therefore we gather that no number which is contained in any of
the following expressions is able to be made into a square.
\[
\begin{array}{rlp{1in}rl}
4mn-&(m+n)&&4mn+&3(m+n)\\
8mn-&(m+n)&&8mn+&7(m+n)\\
8mn-&3(m+n)&&8mn+&5(m+n)\\
12mn-&(m+n)&&12mn+&11(m+n)\\
12mn-&7(m+n)&&12mn+&5(m+n)\\
20mn-&(m+n)&&20mn+&19(m+n)\\
20mn-&3(m+n)&&20mn+&17(m+n)\\
20mn-&7(m+n)&&20mn+&13(m+n)\\
20mn-&9(m+n)&&20mn+&11(m+n)\\
24mn-&(m+n)&&24mn+&23(m+n)\\
24mn-&5(m+n)&&24mn+&19(m+n)\\
24mn-&7(m+n)&&24mn+&17(m+n)\\
24mn-&11(m+n)&&24mn+&13(m+n)\\
28mn-&(m+n)&&28mn+&27(m+n)\\
28mn-&9(m+n)&&28mn+&19(m+n)\\
28mn-&11(m+n)&&28mn+&17(m+n)\\
28mn-&15(m+n)&&28mn+&13(m+n)\\
28mn-&23(m+n)&&28mn+&5(m+n)\\
28mn-&25(m+n)&&28mn+&3(m+n)\\
&&\textrm{etc.}&&
\end{array}
\]
It is to be noted also that in the formulas in the second column
the numbers $m$ and $n$ for the coefficient of $m+n$ must be prime.
This condition requires the restriction that we established in the beginning,
that in the form $aa+Nbb$ the numbers $a$ and $b$ should be relatively prime
numbers: for unless this condition is observed, any number whatsoever would be
able
to be a divisor of this form. As well, it is clear from the preceding that
with this condition observed,
if $4Nmn-A(m+n)$ is not able to be a square, then likewise it is wide open
that $4Nmn-A(m+n) \pm 4Np(m+n)$ is not able to be a square.

\begin{center}{\Large Scholion 3.}\end{center}
Now the expression $aa-Nbb$ shall be considered when it has no divisor
contained in the formula $4Nm\pm A$. It will therefore be $aa-Nbb \neq 4Nmu \pm
Au$ or $aa \neq 4Nmu +NAA \pm Au$. It may be put $NA \pm u=d$, or
$u=\pm d \pm NA$, and it will be
$aa \neq \pm 4Nmd+4NNAm+Ad$; it becomes $d= \pm 4NNn$ and it will become
$16N^3mn \mp 4NNAm \pm 4NNAn \neq aa$, from which it follows that no number
contained in the formula $4Nmn \pm A(m-n)$ can be a square. Neither therefore
will any number contained in the expression $4Nmn \pm A(m-n) \pm 4Np(m-n)$
be able to be a square, but only if the condition recalled from earlier is
observed, that $a$ and $b$ are relatively prime numbers. Consequently
from the preceding theorems the following formulas are deduced, which are never
able to permit square numbers.
\[
\begin{array}{rlp{1in}rl}
8mn\pm&3(m-n)&&8mn\pm&5(m-n)\\
12mn\pm&5(m-n)&&12mn\pm&7(m-n)\\
20mn\pm&3(m-n)&&20mn\pm&17(m-n)\\
20mn\pm&7(m-n)&&20mn\pm&13(m-n)\\
24mn\pm&7(m-n)&&24mn\pm&17(m-n)\\
24mn\pm&11(m-n)&&24mn\pm&13(m-n)\\
28mn\pm&5(m-n)&&28mn\pm&23(m-n)\\
28mn\pm&11(m-n)&&28mn\pm&17(m-n)\\
28mn\pm&8(m-n)&&28mn\pm&15(m-n)\\
&&\textrm{etc.}&&
\end{array}
\]
by attending to which moreover it is easily seen that both of the numbers
$m$ and $n$ must be prime to the coefficient of $(m-n)$: for otherwise,
if the letters given in the formula $12mn \pm 5(m-n)$ were put $m=5p$
and $n=5q$, it would follow that $12\cdot 25pq \pm 25(p-q)$ and thus the
formula $12pq \pm (p-q)$ ought to be a square, which
however is false.
\end{document}